%
%
%
\input amstex
\documentstyle{amsppt}
\def\curraddr#1\endcurraddr{\address
   {\it Current address \/}: #1\endaddress}
 
\topmatter
\title Rad\'o theorem and its generalization\\
for $CR$-mappings
\endtitle
\rightheadtext{RAD\'O THEOREM AND ITS GENERALIZATION}
\author E.M.Chirka \endauthor
\address Steklov Mathematical Institute, Vavilov st. 42, Moscow
GSP-1, 117966 Russia \endaddress
\email chirka\@mph.mian.su \endemail
\abstract  A generalization of Rad\'o's theorem for
CR-functions on locally Lipschitz hypersurfaces is obtained. It
is proved also that closed preimages of pluripolar sets by CR-mappings
are removable for bounded CR-functions.
\endabstract
\thanks
The final version of this paper will be
submitted for publication in "Mat.Sbornik"
\endthanks
\endtopmatter
 
\document
\head 1. Introduction
\endhead
A well-known theorem of Rad\'o \cite{9} states that a continuous
function $f$ defined on a domain in ${\Bbb C}$ and holomorphic on the
complement of its zero set $f^{-1}(0)$ is holomorphic everywhere.
The result is correct for holomorphic functions in ${\Bbb C}^n$
as well as in the plane. It is well-known that
$f^{-1}(0)$ can be replaced by $f^{-1}(E), ~~E$ a closed
subset of zero capacity in ${\Bbb C}$ (see \cite{14}).
 
Recently J.-P.Rosay and E.L.Stout \cite{10} have
shown that an analogue of the classical Rad\'o's theorem take place for
$CR$-functions on a $C^2$-hypersurface in ${\Bbb C}^n$ with nonvanishing
Levi form. Then H.Alexander \cite{1} has proved the removability
in the same situation of closed sets of the type $f^{-1}(E),~E$
a closed polar set in ${\Bbb C}$. We improve here these results in the
following theorem which can be considered as an extension of
Rad\'o's theorem to bounded $CR$-mappings of hypersurfaces.
 
\proclaim{Theorem 1} Let $\Gamma $ be a locally Lipschitz
hypersurface in ${\Bbb C}^{n}$ with one-sided extension property at
each point, $\Sigma$  is a closed subset of $\Gamma$ and
$$f: \Gamma \setminus \Sigma \longrightarrow {\Bbb C}^{m} \setminus E$$
is a $CR$-mapping of class $L^{\infty}$ such that the cluster set of
$f$ on $\Sigma$ along of Lebesque points of $f$ is contained in a closed
complete pluripolar set $E$. Then there is a CR-mapping
$\tilde{f} : \Gamma \longrightarrow {\Bbb C}^{m}$
of class $L^{\infty}(\Gamma)$ such that
$\tilde{f} \mid_{\Gamma \setminus \Sigma} = f$.
\endproclaim
 
We say that $\Gamma$ has one-sided extension property at its point
$a$ if for an arbitrary neighbourhood $U \ni a$ there is a (smaller)
neighbourhood $V \ni a$ and a connected component
$W$ of $V \setminus \Gamma$ such that $a \in \bar{W}$ and every
bounded $CR$-function on $\Gamma \cap U$ extends holomorphically into $W$.
As it was shown by Tr\'epreau \cite{15} this property at each point has
an arbitrary locally Lipschitz hypersurfases in ${\Bbb C}^{2}$ which
contains no analytic discs. The same is true for hypersurfases of
class $C^{2}$ in ${\Bbb C}^{n}$ containing no complex hypersurfaces
(see \cite{15}).
 
We show indead that the trivial extension of $f$ by a constant on
$\Sigma$ is a $CR$-mapping on the whole $\Gamma$.
 
Theorem 1 is true also without the supposition of the one-sided
extension property, but the general case is more compicated, and
this is related with analytic discs belonging to $\Gamma$. We
shall prove the general theorem in the next paper.
 
Theorem 1 is equivalent to the following new result on the removability
of singularities for bounded $CR$-functions.
 
\proclaim{Theorem 2} Let $\Gamma , \Sigma , E$ and $f$ be as in Theorem 1.
Then each $CR$-function of class $L^{\infty}$
on $\Gamma \setminus \Sigma$ extends to a $CR$-function of class
$L^{\infty }$ on $\Gamma$.
\endproclaim
 
Theorem 1 follows obviously from Theorem 2. In opposite direction,
given a $CR$-function $g$ of class $L^{\infty}$
on $\Gamma \setminus \Sigma$ corresponds the mapping
$$(f,g) : \Gamma \setminus \Sigma \longrightarrow
{\Bbb C}^{m+1} \setminus E \times {\Bbb C}$$
which cluster set on $\Sigma$ is contained in the closed complete
pluripolar set $E \times {\Bbb C}$. By Theorem 1, the map $(f,g)$
extends to a $CR$-map of whole $\Gamma$, and thus its last
component $g$ extends to there as a $CR$-function.
 
If $\Gamma$ is the boundary of a bounded domain or if $\Gamma$ admits
one-sided holomorphic extension of $CR$-functions (say, if
$\Gamma \in C^{2}$ contains no complex hypersurface, as in
\cite{15}) then it follows from Theorem 2 and a uniqueness theorem that
the Hausdorff $(2n-1)$-measure of $\Sigma$ vanishes, and
$\Gamma \setminus \Sigma$ is locally connected.
 
We prefer to work here with the class $L^{\infty}$ instead of $C$,
since $L^{\infty}$ is stable by considering extensions whereas the
$CR$-extension of a bounded continuous $CR$-function from
$\Gamma \setminus \Sigma$ onto $\Gamma$ is not continuous in general.
 
The proof of Rad\'o theorem for $CR$-functions in \cite{10} is based on
results \cite{8} on the holomorphic continuation of
$CR$-functions from a part of the boundary of a domain in
${\Bbb C}^n,~ n \ge 2$ (see also \cite{7}). In the proof
of Theorem 1 we use instead of this a geometric extension of the graph
in spirit of R.Harvey and H.B.Lawson \cite{6}.
Our starting point was a generalization of the Harvey~-~Lawson theorem
\cite{6} on boundaries of holomorphic chaines for MC-cicles in the
complement to a polynomially convex compact set \cite{4\rm, 1985}.
It can be considered as a geometric version (for $CR$-functions of class
$C^1$) of theorems on holomorphic continuation in
\cite{8},\cite{10},\cite{7}.
 
Let us specify the terminology.
 
We say that a hypersurfase $M$ in a smooth $k$-dimensional manifold
$\Cal M$ is {\sl locally Lipschitz} if for every point $a \in M$ there
is a coordinate chart $(U,x), ~x = (x', x_{k})$ on $\Cal M$ such that
$a \in U$ and $M \cap U$ is represented as the graph $x_{k} = h(x')$
of a function $h$ over the domain in ${\Bbb R}^{k-1}$ which sutisfies
there Lipschitz condition $\mid h(b) - h(c) \mid \le C \mid b - c \mid$
with a constant $C$.  Note that the Hausdorff $k-1$-measure
(with respect to some fixed smooth metric on $\Cal M$) restricted
to such $M$ is locally finite, and $M$ has tangent planes in almost
every point with respect to this measure (Rademacher's theorem, see
e.g. \cite{5\rm,3.1.6.}). Thus the integral on $M$ for a differential
$(k-1)$-form $\varphi$ with Lipschitz coefficients and with compact
$supp{\varphi} \cap M$ is well-defined.
 
A point $a$ in a locally Lipschitz $M \subset \Cal {M}$ is called a
{\sl Lebesque point} for a given vector-function $f$ of class
$L^{1}_{loc}(M)$ with values in ${\Bbb R}^N$ if there is a constant
$\tilde{f}(a) \in {\Bbb R}^N$ such that
$$
r^{-k+1}\int_{M \cap \mid x' \mid < r}
\mid f(x) - \tilde{f}(a) \mid \,dx' \longrightarrow 0
$$
as $r \rightarrow 0$ in the chart $(U,x)$ with $x(a) = 0$ described
above. It is wellknown (see \cite{5\rm,2.9.8.}) that almost every
point $a \in M$ is such a point and $f(a) = \tilde{f}(a)$ almost everywhere.
Thus we shall assume in the further that $f$ is defined (as $\tilde{f}(a)$)
on the set of its Lebesque points $a$ {\sl only}.
 
For a locally Lipschitz hypersurfase $M$ in a complex $n$-dimensional
manifold $\Cal M$ the notion of $CR$-functions of class $L^{1}_{loc}(M)$
is well-defined:
a function $f$ of this class is a $CR$-function on $M$ if
$\int_{M} f \bar {\partial } \varphi = 0$
for every smooth form $\varphi $ of bidegree $(n,n-2)$ in $\Cal M$
with compact $supp \varphi \cap M$.
 
A set $E \subset {\Bbb C}^{m}$ is called complete pluripolar, if there
is a plurisubharmonic function $\varphi$ in ${\Bbb C}^{m}$ such that
$E = \{\zeta : \varphi(\zeta) = - \infty \}$.
\head 2. One-sided holomorphic extension
\endhead
The problem is local, so we can assume that $\Gamma \ni 0$ is represented
as the graph $v = h(z^{'},u)$ of some Lipschitz function in a domain
in the space of variables $(z_{1},...,z_{n-1},Re z_{n}) =  (z^{'},u)$
(it is convenient to use the notation $z_{n} = u + iv$).
 
Fix a connected component $\Gamma_{0}$ of $\Gamma \setminus \Sigma$,
set $f_{1} = f$ on $\Gamma_{0}, f_{1} = 0$ on
$\Gamma \setminus \Gamma_{0}$ and denote by $\Gamma_{1}
\supset \Gamma_{0}$ the set of points $a \in \Gamma$ such that
$f_{1}$ is a $CR$-function in a neighbourhood of $a$ on $\Gamma$.
We have to show that $\Gamma_{1} = \Gamma$.
 
By the one-sided extension
property, for each point $a \in \Gamma_{1}$ there is a neighbourhood
$V_{a} \ni a$ and a connected component $W_{a}$ of
$V_{a} \setminus \Gamma$ such that $a \in \overline{W_a}$ and each
$CR$-function on $\Gamma_{1}$
extends holomorphically into $W_a$. Shrinking $V_a$ we can
assume that the intersection of $V_a$ with each line
$(z^{'},u) = const$ is an interval
(i.e. $V_a$ is convex in $v$-direction) intersecting
$\Gamma_{1}$. Then the union of all $W_{a}, a \in \Gamma_{1}$ is an
open set of the form $W^{+} \cup W^{-}$ where $W^{+}$ is plased over
$\Gamma$ (i.e. $v > h(z^{'},u)$ on $W^{+}$) and $W^{-}$ is contained in
$\{v < h(z^{'},u)\}$. It follows that the set
$$W = W^{+} \cup W^{-} \cup (\overline{W^{+}} \cap \overline{W^{-}} \cap
\Gamma_{1})$$
is open, convex in $v$-direction, and $\bar{W} \supset \Gamma_{1}$.
By a uniqueness theorem and a removable singularities theorem the
holomorphic extensions of $f_{1}$ into $W_{a}, a \in \Gamma_{1},$
constitute holomorphic functions in $W^{+}$ and $W^{-}$, and these
functions extend to a holomorphic (vector-) function in $W$ which
we denote by the same symbol $f_{1}$.
 
By the construction, there is a Lipschitz function
$\epsilon(z^{'},u)$ such that $\epsilon = 0$
outside of $(z^{'},u)(\Gamma_{1})$, the hypersurface
$\Gamma^{'} : v = (h + \epsilon)(z^{'},u)$ is contained in
$W \cup \Sigma_{1}$, and $\Gamma^{'} \setminus
\Sigma_{1}$ is a smooth ($C^{\infty}$ or even $C^{\omega}$ if you
want). The set $(f_{1})^{-1}(E) \cap W$ is pluripolar, so we can
assume that its intersection with $\Gamma^{'} \setminus
\Sigma_{1}$ has the Hausdorff dimension $2n - 3$, in particular,
it has the locally connected complement in
$\Gamma_{'} \setminus \Sigma$. Set $\Sigma_{1}^{'} = \Sigma_{1} \cup
(\Gamma_{1}^{'} \cap (f_{1})^{-1}(E)$
and $\Gamma_{1}^{'} = \Gamma^{'} \setminus \Sigma_{1}^{'}$.
If $\epsilon(z^{'},u)$ is taken sufficiantly small and rapidly tends
to zero as $(z^{'},u)$ approaches to $(z^{'},u)(\Sigma_{1})$, then
$$
f_{1} \mid \Gamma_{1}^{'} \longrightarrow {\Bbb C}^{m} \setminus E
$$
and the cluster set of $f_{1}$ on $\Sigma_{1}^{'}$ is contained
in $E$. Thus, substituting $\Gamma$ onto $\Gamma^{'}, \Sigma_{1}$
onto $\Sigma_{1}^{'}$ and $\Gamma_{1}$ onto $\Gamma_{1}^{'}$ and
then restoring old notations, we can assume that $\Gamma_{1}$
is smooth and the mapping $f_{1}$ is holomorphic in a
neighbourhood of $\Gamma_{1}$.
\head 3. Reducing to $n = 2$
\endhead
We assume as above that $\Gamma \ni 0$ is represented
as the graph of some Lipschitz function over a domain in the space
of variables $z_{1},...,z_{n-1},Re z_{n}$. Then the vector $(0,...,0,i)$
does not belong to $C_{a}\Gamma$, the tangent cone to $\Gamma$ at the
point $a$, for all $a \in \Gamma$. Shrinking $\Gamma$ a little we can
assume that the same is true for some ${\Bbb C}$-linearly independent
system of vectors $\xi_{1},...,\xi_{n}, ~~\xi_{j} \not\in C_{a}\Gamma$
for $a \in \Gamma, j = 1,...,n$. Making a suitable ${\Bbb C}$-linear
changing of coordinates we obtain the situation when
$ie_{j} \not\in C_{a}\Gamma$ for all standart coordinate orts
$e_{j}$ in ${\Bbb C}^n$.
It follows then that for each $j$ there is a neighbourhood $U_{j} \ni 0$
such that $\Gamma \cap U_{j}$ is represented as the graph of a Lipschitz
function over a domain in the space of variables
$z_{k}, k \ne j, ~ Re z_{j}$. Set $U = \cap_1^n U_{j}$.
 
We have to show that $\int_{\Gamma} f_{0} \bar{\partial} \varphi = 0$ for
an arbitrary smooth $(n,n-2)$-form $\varphi$ with $supp \varphi
\subset U$. This form is represented as $\sum_{j < k} \varphi_{jk}
dz_{j} \wedge dz{k} \wedge dV_{jk}$ where $\varphi_{jk}$ are smooth
functions supported in $U$ and
$dV_{jk} = \prod_{l \ne j,k} i dz_{l} \wedge d \bar{z}_{l}$.
 
By the construction, the projection $\Gamma_{jk}$ of $\Gamma \cap U_{k}$
into the space of variables $\{z_{l}, l \ne j,k\}$ is an open set, and
$\Gamma_{c(j,k)} = \Gamma \cap \{ z_{l}=c_{l},~ l \ne j,k \}$ is a
Lipschitz hypersurface in $\{ z_{l}=c_{l},~ l \ne j,k \} \simeq {\Bbb C}^2$
for all $c(j,k) \in \Gamma_{jk}$. As
$\bar{\partial} \varphi = \sum_{j < k} \bar{\partial} (\varphi _{jk}
dz_{j} \wedge dz{k}) \wedge dV_{jk}$, we have by Fubini theorem for
differential forms (see e.g. \cite{4\rm, A4.4.}) that
$$\int_{\Gamma} f_{0} \bar{\partial} \varphi =
\sum_{j < k} \int_{\Gamma_{jk}} (\int_{\Gamma_{z(j,k)}} f_{0}
\bar{\partial} (\varphi _{jk} dz_{j} \wedge dz{k})) \wedge dV_{jk}.$$
If $f_{0} \mid \Gamma_{c(j,k)} \in CR(\Gamma_{c(j,k)})$ for almost every
$c(j,k) \in \Gamma_{jk}$,then almost all inner integrals vanish, and
the righthand side is zero.
 
Taking $\varphi$ in a dense sequence of such forms with compact
$supp \varphi \cap \Gamma \setminus \Sigma$ we obtain from this
representation that $$f \mid \Gamma _{c(j,k)} \setminus \Sigma
\longrightarrow {\Bbb C}^{m} \setminus E$$
are the mappings of the class $CR \cap L^{\infty}$ for almost every
$c(j,k) \in \Gamma_{jk}$. As
$\Gamma_{jk} \subset \{ z_{l} = c_{l} , l \ne j,k \} \simeq {\Bbb C}^2$,
we obtain that it is enough to prove Theorem 1 for the case $n = 2$.
\head 4. Analytic extension of the graph
\endhead
To show that $\Sigma_{1}$ is empty we assume that $0$ is the
boundary point of $\Gamma_{1}$ in $\Gamma$ and come at last
to a contradiction.
 
The base domain $G \subset {\Bbb C} \times \Bbb{R}$  can be taken
bounded and convex, and the function $h(z,u)$ defined and with
Lipschitz condition in a neighbourhood of $\bar{G}$. Then the
graph $S : v = h(z,u)$ over $bG$ is a two-dimensional sphere in
${\Bbb C}^{2}$. As $0$ is limiting point for $\Gamma_{1}$, we can
assume, that $S$ is not contained in $\Sigma_{1}$. By Shcherbina's
theorem \cite{12, 13} the polynomially convex hull $\tilde{S}$
of $S$ is the graph of a continuous function $\tilde{h}(z,u)$ over
$\bar{G}$ foliated  in a one-parametric family of analytic discs
with boundaries on $S$.
 
The graph $M$ of the map $f_{1}$ over $\Gamma_{1}$ is a smooth
maximally convex 3-dimensional manifold in ${\Bbb C}^{2} \times
{\Bbb C}^{m}$ which boundary $\bar{M} \setminus M$ is contained in
$(\tilde{S} \times {\Bbb C}^{m}) \cup ({\Bbb C}^{2} \times E)$. As
$f_{1}$ is uniformly bounded, there are closed balls $B_{2},
B_{m}$ with centers in origins such that $\bar{M} \setminus M$
is contained in $(\tilde{S} \times B_{m}) \cup (B_{2} \times E)$.
This compact set is polynomially convex due to the following.
 
\proclaim{Lemma 1} Let $X_{1} \subset X_{2}$ be polynomially
convex compact sets and $Y$ is a complete pluripolar set in
${\Bbb C}^{N}$. Then the set $X = X_{1} \cup (Y \cap X_{2})$ is
polynomially convex.
\endproclaim
 
\demo{Proof} The set $Y$ is represented as $\{\zeta :
\varphi (\zeta) = - \infty \}$ for some function $\varphi$
plurisubharmonic in ${\Bbb C}^{N}$. If $a \not\in X_{1} \cup Y$
then there is a polynomial $p$ such that $p(a) = 1$ and
$\mid p \mid < 1$ on $X_{1}$. Let $C = sup\{\varphi (\zeta) :
\zeta \in X_{1} \}$ and a positive integer $s$ is taken so big
that $\mid p(\zeta) \mid ^{s} e^{C} < e^{\varphi (a)}$ for all
$\zeta \in X_{1}$ (it is possible because $a \ni Y$). Then the
function $\psi = \mid p \mid ^{s} e^{\varphi}$ is
plurisubharmonic in ${\Bbb C}^{N}, \psi (\zeta) < \psi (a)$
for $\zeta \in X$, and the same is true for $\zeta \in Y$
because $\psi \mid Y = 0$. It follows from the maximum
principle for plurisubharmonic functions on polynomially
convex hulls (see, e.g. \cite{3}) that $a$ is not contained
in the hull of $X$. The rest follows from the inclusion
$X \subset X_{2}$.
\enddemo
 
Thus, the polynomially convex hull of the set
$(S \times B_{m}) \cup (B_{2} \times E)$ for $B_{2} \supset S$
is the compact set
$$K = (\tilde{S} \times B_{m}) \cup (B_{2} \times E),$$
and the graph $M$ of $f_{1}$ is attached to this $K$. By a
generalization of Harvey~-~Lawson theorem in \cite{4\rm,
Theorem 19.6.2} there is a two-dimensional (complex) analytic
subset $A$ in ${\Bbb C}^{2+m} \setminus (K \cup M)$ such that
$A \cup K \cup M$ is compact and $M \setminus K \subset \bar{A}$.
 
\head 5. The projection of the extension
\endhead
We show that $A$ is the graph of a holomorphic mapping over an
open set in ${\Bbb C}^{2}$ with the boundary in $\Gamma \cup \tilde{S}$.
(The main difficulty here is that the projection of $\bar{A}$ is
as well as $\Gamma \cup \tilde{S}$ contained in the ball $B_{2}$,
the "shadow" of $B_{2} \times E$.) It is convenient to use in
this Section coordinates $(z^{'},z^{"})$ for ${\Bbb C}^{2} \times
{\Bbb C}^{m}$.
 
We essentially use the pluripolarity of $B_{2} \times E$ and the
following result due to E.Bishop \cite{2, 11} on the removability
of pluripolar singularities for analytic sets (see \cite{4}).
 
\proclaim{Lemma 2} Let $Y$ be a closed complete pluripolar subset
of a bounded domain $U = U^{'} \times U^{"} \subset {\Bbb C}^{n+m}$,
and $A$ is a pure p-dimensional analytic subset in $U \setminus Y$
without limit points on $U^{'} \times bU^{"}$. Suppose that $U^{'}$
contains a nonempty subdomain $V^{'}$ such that $\bar{A} \cap
(V^{'} \times U^{"})$ is an analytic set. Then $\bar{A} \cap U$
is analytic in $U$.
\endproclaim
 
First of all, we apply this lemma to the unbounded component
$U^{'}$ of ${\Bbb C}^{2} \setminus (\tilde{S} \cup \Gamma)$. Let
$U^{"}$ be an open ball in ${\Bbb C}^{m}$ containing $B_{m}$ and
$Y = U \cap B_{2} \times E$. Then $A \cap (U \setminus Y)$ is
an analytic set sutisfying the conditions of Lemma 2. By the
maximum principle, $A$ is projected into $B_{2}$ because its
boundary is contained in $K \cup M$. Thus, for $V^{'} = U^{'}
\setminus B_{2}$, the set $\bar{A} \cap (V^{'} \times U^{"})$
is empty (hense analytic). It follows from Lemma 2 that
$\bar{A} \cap U$ is analytic in $U$. As
$\bar{A} \cap (V^{'} \times U^{"})$ is empty and the projection
of $\bar{A} \cap U$ into $U^{'}$ is proper, the set
$\bar{A} \cap U$ is also empty. Thus, we have proved that the
projection of $A$ into ${\Bbb C}^{2}$ is contained in the closure
of the union of all bounded component of ${\Bbb C}^{2} \setminus
(\Gamma \cup \tilde{S})$.
 
Take now an arbitrary point $a^{'} \in \Gamma_{1} \setminus \tilde{S}$
and show that the set $A \cap \{z^{'} = a^{'}\}$ is empty. This
set is closed analytic in $a^{'} \times ({\Bbb C}^{m} \setminus
(E \cup a^{"}))$ where $a^{"} = f_{1}(a^{'})$. As $E \cup a^{"}$
is complete pluripolar, its intersection with $B_{m}$ is
polynomially convex. As $\bar{A} \cap \{z^{'} = a^{'}\}$ is
compact, it follows from a maximum principle on analytic sets
(see, e.g., \cite{4\rm, 6.3.}) that the dimension of
$A \cap \{z^{'} = a^{'}\}$ is zero, i.e. this set is discrete.
Thus, given $b =(a^{'},b^{"}) \in A$ there is a neighbourhood
$U = U^{'} \times U^{"}$ such that the projection of
$A \cap U$ into $U^{'}$ is an analytic covering (see \cite{4}).
But $dim A = 2$, and there is no point in $A \cap U$ over
unbounded componenet of ${\Bbb C}^{2} \setminus (\Gamma \cup \tilde{S})$
which has nonempty intersection with $U^{'}$. This contradiction
shoes that there is no such points $b$, i.e. $A \cap U$ is empty.
 
Let now $U^{'}$ be a bounded component of
${\Bbb C}^{2} \setminus (\Gamma \cup \tilde{S})$ such that
$bU^{'} \cap (\Gamma_{1} \setminus \tilde{S}$ is not empty, and
$a^{'}$ is a point in this nonempty set. Then (see Sect.1.) there
is a neighbourhood
$V^{'} \ni a^{'}$ such that $f_{1}$ is holomorphic in $V^{'}$.
We have in $(V^{'} \times {\Bbb C}^{m}) \setminus M$ two analytic
sets, $A \cap (V^{'} \times {\Bbb C}^{m})$ and the graph of $f_{1}$
over $V^{'} \cap U^{'}$, of pure dimension $2$ with the same
smooth boundary $M \cap (V^{'} \times {\Bbb C}^{m})$. By a boundary
uniqueness theorem for analytic sets (see \cite{4\rm, 19.2.}),
these sets coinside. Thus, the analytic covering
$A \cap (U^{'} \times {\Bbb C}^{m}) \longrightarrow U^{'}$ is
onesheeted over $V^{'} \cap U^{'}$, which follows that it is
one-sheeted over whole $U^{'}$. It means that $A$ over $U^{'}$
is the graph of a bounded holomorphic map, and this map is a
continuation of $f_{1}$ into $U^{'}$. In particular, we obtain
that $f_{1}$ as the boundary value of this map is $CR$ on
$bU^{'} \cap (\Gamma \setminus \tilde{S})$.
 
In terms of components of $\Gamma \setminus \tilde{S}$, it means
that there are only two possibilities: either this component is
contained in $\Gamma_{1}$ or it is contained in $\Sigma_{1}$.
 
\head 6. Removability of $\Sigma$
\endhead
Return to notations $(z,w)$ for coordinates in ${\Bbb C}^{2}$ and
denote by $p$ the projection $(z,w) \mapsto (z,u)$ into
${\Bbb C}^{2} \times \Bbb{R}$.
 
Let $\delta \subset \tilde{S}$ be an analytic disc with the boundary
in $S$ such that $p(\delta) \cap p(\Gamma_{1})$ is not empty.
Show that $p(\delta) \subset p(\Gamma_{1})$.
 
Suppose it is not. Then there is a point $a \in \Sigma_{1}$
such that $p(a) \in p(\delta)$, and a convex domain
$G_{1} \subset G$ containing $p(a)$ such that
$bG_{1} \cap p(\delta) \cap p(\Gamma_{1})$ is not empty, say,
it contains $p(b)$ for some $b \in \Gamma_{1}$. By Sect.1, $f_{1}$
is holomorphic in a one-sided neighbourhood $V$ of $b$ convex
in $v$-direction, as in Sect.1. Let $h_{1}$ be a Lipschitz
function on $bG_{1}$ such that its graph $S_{1} : v = h_{1}(z,u)$
is contained in $\Gamma \cup V$ but does not contain $b$. We can
assume $v < h(z,u)$ on $S_{1} \cap V$ (changing $w$ onto $- w$
and shrinking $V$ if it is nessessary ). Then for
$h_{t} = th + (1-t)h_{1}, 0 < t < 1,$ we have $h_{t} \leq h$ and
$h_{t} < h$ over $p(V)$. By Shcherbina's theorem \cite{12, 13}
polynomially convex hull $\widetilde{S_{t}}$ of
$S_{t} : v = h_{t}(z,u)$ is the graph of a continuous function
$\tilde{h}_{t}$ over $G_{1}$ foliated in a one-parametric family of
analytic discs with boundaries in $S_{t}$. As $h_{t} \leq h$,
we have $\tilde{h}_{t} \leq \tilde{h}$, and
$\tilde{h}_{t} \rightarrow \tilde{h}$ as $t \rightarrow 0$ by continuity.
 
Thus, for $t > 0$ small enough there is a disc
$\delta_{t} \subset \widetilde{S_{t}}$ such that
$p(a) \in \delta_{t}$ and $p(b\delta) \cap p(V)$ is not empty.
As $v \leq h(z,u)$ on $\delta_{t}$ and $v < h(z,u)$ on
$b\delta \cap V$, we have $v < \tilde{h}(z,u)$ on the whole
$\delta_{t}$, i.e. the disc $\delta_{t}$ is placed strongly
under the hypersurfase $\tilde{S}$. (This is true because the discs
$\delta_{t} - (0, i\epsilon)$ for $\epsilon > 0$ do not
intersect $\tilde{S}$, and we can apply the argument principle.)
 
As we proved above (with $\widetilde{S_{t}}$ instead of $\tilde{S}$),
each component of $p(\delta_{t} \setminus \Gamma)$
is either contained in $p(\Gamma_{1})$ or it is contained in
$p(\Sigma_{1})$. But $p(\delta_{t} \cap \Gamma)$ and
$p(\tilde{S} \cap \Gamma)$ have no common point because
$\delta_{t} \cap \tilde{S}$ is empty, and the projection
$p \mid \Gamma$ is one-to-one. As
$p(\delta_{t}) \cap p(\Gamma_{1})$ is not empty by the
construction, the set $p(\delta_{t}) \cap p(\Sigma_{1})$ must to
be empty, in particular, $p(a) \not\in p(\Sigma_{1})$. The
contradiction (with choosing of $a$) shows that there is no such
point $a$, i.e. $p(\delta) \subset p(\Gamma_{1}$.
 
If $\Sigma_{1}$ is not empty, it follows from the above that
there is a disc $\delta^{0} \subset \tilde{S}$ such that
$p(\delta^{0}) \subset p(\Sigma_{1}) \cap \partial(p(\Gamma_{1}))$.
As $\delta^{0}$ is not contained in $\Gamma$, there is a point
$c \in \Sigma_{1} \setminus \tilde{S}$ such that
$p(c) \in p(\delta^{0})$. The component $\Gamma_{c}$ of
$\Gamma \setminus \tilde{S}$ containing $c$ has nonempty intersection
with $\Gamma_{1}$ because $c$ is limiting poimt for $\Gamma_{1}$.
As it was proved above, it follows that $\Gamma_{c}$ is contained
in $\Gamma_{1}$. This contradiction (with $c \in \Gamma_{1}$)
shows that $\Sigma_{1}$ is empty.
 
Thus, we have proved that for each component $\Gamma_{0}$ of
$\Gamma \setminus \Sigma$ the map $f_{1}$ (equal to $f$ on
$\Gamma_{0}$ and to $0$ outside of it) is $CR$ on the whole
$\Gamma$. As $\Gamma$ contains no analytic discs, this $f_{1}$
extends holomorphically into one-sided neighbourhouds of each
point of $\Gamma$. It follows from a boundary uniqueness theorem
for holomorphic functions that $\Gamma \setminus \Gamma_{0}$
has zero Hausdorff $3$-measure, in particular, $\Gamma_{0}$ is
the single component of $\Gamma \setminus \Sigma$. In other
words, the set $\Sigma$ has zero Hausdorff $(2n-1)$-measure
for general $n$, and its complement in $\Gamma$ is locally
connected (for every connected open set
$\Gamma^{'} \subset \Gamma$ the set $\Gamma^{'} \setminus \Sigma$
is also connected).
\end